\documentclass[12pt]{iopart}

\usepackage{iopams}  

\usepackage{graphicx}

\usepackage{amssymb,amsfonts}

\newtheorem{theorem}{Theorem}
\newtheorem{proposition}{Proposition}
\newtheorem{lemma}{Lemma}
\newtheorem{corollary}{Corollary}
\newtheorem{remark}[corollary]{Remark}

\newcommand{\dist}{\mbox{\rm dist}}

\newcommand{\QED}{\hfill \raisebox{-2pt}{\rule{5.6pt}{8pt}\rule{4pt}{0pt}}%
  \smallskip\par}

\begin{document}

\title{An iterative algorithm for the square-root Lasso}

\author{Patrizia Boccacci $^{(*)}$ \& Christine De Mol $^{(+)}$ \& Ignace Loris $^{(+)}$}

\address{ $^{(*)}$ DIBRIS, Universit\`a di Genova, Italy, and $^{(+)}$ Department of Mathematics, Universit\'e libre de Bruxelles, Belgium}
\ead{patrizia.boccacci55@gmail.com, Christine.De.Mol@ulb.be, Ignace.Loris@ulb.be}
\vspace{10pt}
\begin{indented}
\item[]October 2025
\end{indented}

\begin{abstract} In the framework of sparsity-enforcing regularisation for linear inverse problems, we consider the minimisation of a square-root Lasso cost function.  To solve this problem we devise a simple modification (called SQRT-ISTA) of the Iterative Soft-Thresholding Algorithm (ISTA) for the Lasso problem and we prove convergence for this algorithm. Under some additional assumptions, we derive an upper bound on the convergence rate of the cost function. We also generalise these results to the case of the group square-root Lasso, where sparsity is enforced for groups of variables instead of individual ones. 

\end{abstract}

\section{Introduction}
\label{intro}

We consider the following linear inverse problem: find an object of interest $f$ from the observation or measurement of $g=Af$ where $A$ is a given linear operator from a space $\mathcal{F}$, the object space, to which $f$ belongs, into a data or image space $\mathcal{G}$, to which $g$ belongs. The spaces $\mathcal{F}$ and $\mathcal{G}$ are assumed to be Hilbert spaces (typically $L^2$ spaces) and $A:\mathcal{F} \to \mathcal{G}$ to be a bounded linear operator. In the special instance of a finite-dimensional problem, $A$ is just a matrix and $\mathcal{F}$ and $\mathcal{G}$ Euclidean spaces. For simplicity, we will only consider real-valued $f$ and $g$.

Since measurements are affected by noise and the corresponding inverse problem is typically ill-posed, it should be regularised. A usual strategy to achieve this is to replace it by a variational problem. In the framework of sparsity-enforcing regularisation, the specific variational problem we consider in the present work is to find a minimiser $f_*$ for the following cost function 

\begin{equation}
 \label{phi}
\Phi_{\mu}(f)
= \| Af-g \| + \mu\, \Vert f\Vert_1
\end{equation}
where $\mu>0$ is a fixed regularisation parameter and $\Vert f\Vert_1$ is a sparsity-enforcing norm defined by
\begin{equation} \label{pen_w}
\Vert f\Vert_1 =  
\sum_{\gamma \in \Gamma} 
 w_\gamma\, \vert f_\gamma|.
\end{equation}
By  $f_\gamma$ we denote the scalar product $f_\gamma = (f,\varphi_\gamma)$ where $\{\varphi_\gamma\}_{\gamma\in \Gamma}$ is a fixed -- but otherwise arbitrary --
orthonormal basis of the object space $\mathcal{F}$. The ``residual'' (or  ``discrepancy'') $\| Af-g \|$ is measured by means of the norm of the image/data Hilbert space 
$\mathcal{G}$. It is well known that if $f$ is expected to be sparse in the given basis, i.e. if many coefficients $f_\gamma$ are expected to be zero, the penalty (\ref{pen_w}), which is a weighted $\ell^1$-norm on the sequence of  coefficients $f_\gamma$ on that basis, will favour the recovery of such sparse objects. For simplicity, we will consider the case of unit weights $w_\gamma=1$ but this is not restrictive; the extension to arbitrary positive weights uniformly bounded from below $(w_\gamma \geq c > 0)$ is straightforward and is needed when dealing with wavelet bases. 

Dating back to \cite{Daubechies2004}, sparsity-enforcing regularisation has been abundantly studied in the literature but mostly resorting as in  \cite{Daubechies2004} to the minimisation of the following cost function
\begin{equation}
 \label{phi_Lasso}
\Phi_{\tilde \mu}^{\rm Lasso}(f)
= \| Af-g \|^2 + \tilde\mu\, \Vert f\Vert_1
\end{equation}
often referred to as the ``Lasso'' cost function. We notice that this latter cost function is not scale invariant, in the sense that if $g$ is rescaled as $sg$ with a  scale parameter $s>0$, as e.g. by a change of measurement unit, the corresponding minimiser of 
(\ref{phi_Lasso}) is not $sf$. This may appear as a drawback of the Lasso cost function, which looks somewhat awkward from a dimensional point of view. The cost function (\ref{phi}), however, is scale invariant as is the usual Tikhonov cost function where the sparsity-enforcing penalty (\ref{pen_w}) in (\ref{phi_Lasso}) is replaced by the quadratic penalty $\Vert f \Vert^2$.  The minimisation of (\ref{phi}) is usually called the ``square-root Lasso''. While its scale invariance appears as an advantage, taking the square-root of the quadratic residual renders it nonsmooth -- as is the $\ell^1$-norm penalty (\ref{pen_w}) -- and hence makes the optimisation problem more difficult. 

We should remark though that both the Lasso and the square-root Lasso are convex problems. Hence they should be equivalent in the sense that a minimiser of (\ref{phi}) for a given value of $\mu$ should be a minimiser of (\ref{phi_Lasso}) for some value of $\tilde\mu$. Also equivalent are the two associated constrained problems: minimise the residual subject to an upper bound on the norm (\ref{pen_w}) or minimise the norm (\ref{pen_w}) subject to an upper bound on the residual. However, the correspondence between the parameters involved in these constrained problems (between upper bounds and regularisation parameters) is not in general explicit and might be hard to compute. 
We will nevertheless show later in the paper that, for a given minimiser $f_*$, the correspondence between the respective square-root Lasso and Lasso parameters $\mu$ and $\tilde\mu$ is $\tilde\mu=2 \mu \| Af_*-g \|$.

As the acronym Lasso (for Least absolute shrinkage and selection operator) created by Tibshirani \cite{Tibshirani1996}, the appellation ``square-root Lasso'' (abbreviated as SQRT-Lasso or SR-Lasso) originated in statistics in a paper by Belloni et al. \cite{Belloni2011}, in the framework of (finite-dimensional) linear regression, where (\ref{phi}) is usually written in the traditional notations of the field as 
\begin{equation}
 \label{phi_reg}
\Phi_{\mu}(\beta)
= \| X\beta-y \| + \mu\, \Vert \beta\Vert_1.
\end{equation}
$X$ denotes the observed design matrix, $y$ the observed response and $\beta$ the unknown regression coefficient. The $\ell^1$-norm penalty on $\beta$, $\Vert \beta\Vert_1 = \sum_j \vert \beta_j \vert$, enforces its sparsity, allowing for variable selection. 
In that seminal paper, the method is presented as a remedy to the fact that consistency results for the Lasso require a choice of the regularisation parameter $\mu$ that depends on an ``oracle'' knowledge of the variance of the noise (variance which is also hard to estimate). Indeed, the authors show that the square-root Lasso yields a ``pivotal'' estimator, in the sense that the same consistency results as for the Lasso can be obtained with a choice of $\mu$ that is independent of the variance of the noise.  

After \cite{Belloni2011} the square-root Lasso has been extensively studied in the (high-dimensional) statistics literature \cite{VandeGeer2016, Stucky2017, Derumigny2018, Tian2018, Giraud2015, Berk2024}, and the results have also been extended to the group square-root Lasso \cite{Bunea2014}, which, similarly to the group Lasso \cite{Yuan2005, Giraud2015, VandeGeer2016, Stucky2017}, uses a penalty which favours the sparsity of groups of variables, i.e. the sum of the Euclidean norms (not squared) of each group. Its advantages have also been investigated in the framework of ``Compressed Sensing'' \cite{Adcock2021, Adcock2024, Mayrink2024, Adcock2025} or of related sparse recovery problems \cite{Adcock2019}.  We will not comment further on the pivotal properties derived in such settings since they require hypotheses, such as e.g. the restricted isometry property (RIP), the null-space property (NSP) or the restricted eigenvalue property (RE), which do not generally hold for the case of inverse imaging problems.  An exception is one of the algorithms analysed in \cite{Mayrink2024} (Theorems 6 and 19) which does not require such hypotheses.  

We believe that the simple algorithm we analyse here in the framework of inverse problems can be of interest not only for imaging problems (such as e.g. neuroimaging \cite{Bertrand2019}) but also for these applications in statistics and compressed sensing.  Because of the equivalence explained above with the Lasso, we do not of course expect an improvement in the quality of the reconstructions obtained by the square-root Lasso compared to the Lasso. However, the choice of the regularisation parameter could be eased and be less dependent on the effective noise level. 

To the best of our knowledge, our algorithm differs from the other algorithmic solutions proposed in the literature: conic programming \cite{Belloni2011}, alternating direction method of multipliers (ADMM)  \cite{Li2015}, coordinate descent \cite{Ndiaye2017, Zhao2021}, proximal gradient descent and proximal Newton method \cite{Li2019}, and iteratively reweighted least squares (IRLS) \cite{Mayrink2024}.
An exception is the paper \cite{Bunea2014} where an algorithm similar to ours is proposed for the group square-root Lasso and convergence of the sequence of iterates is claimed. However, unfortunately, the proof of convergence provided in the Appendix of that paper appears to be incorrect for the reasons that will be explained in the following. 

In Section \ref{algorithm}, we introduce the proposed iterative algorithm, establishing a connection with the well-known ``Iterative Soft-Thresholding Algorithm'' (ISTA) for Lasso.  Then we study its convergence properties and prove, in finite-dimensional situations, that the sequence of the iterates always converges and that it is guaranteed to converge to a minimiser of the cost function (\ref{phi}) provided the limit residual differs from zero (this is the usual case of interest in inverse problems). We also derive an upper bound on the difference between the cost function and its limit value, as a function of the number of iterations. In Section \ref{group_lasso}, we generalise the algorithm to the case of group square-root Lasso
and show that the proof of convergence of Section \ref{algorithm} extends to this case in a straightforward way. Finally, Section \ref{conclusion} contains some concluding remarks.

\section{A simple iterative algorithm for the square-root Lasso}
\label{algorithm}
\subsection{Derivation of the algorithm}
Notice that 
\begin{equation}
\label{var}
{\rm arg min}_f \Phi_{\mu}(f)
=  {\rm arg min}_{f ,\sigma > 0}  \Phi_{\mu}(f,\sigma)
\end{equation}
where
\begin{equation}
 \label{phisig}
 \Phi_{\mu}(f,\sigma)
=   \frac{\| Af-g \|^2}{\sigma} + \sigma + 2\mu\, \Vert f\Vert_1 
\end{equation}
provided that the minimum is attained for a non-zero value of $\sigma$.
 Inspired by ideas in robust statistics \cite{Huber1981}, this trick has been put forward in \cite{Antoniadis2010} and used for the square-root Lasso in \cite{Giraud2015, VandeGeer2016}.  Notice that the cost function (\ref{phisig}) is jointly convex in $f$
and $\sigma > 0$ (see e.g.  \cite{Boyd2004}, Example 3.18).

This suggests the following alternating minimisation, starting from an arbitrary initial guess for $f_0$, and the corresponding residual 
$\sigma_0= \Vert Af_0-g\Vert$ (e.g. $f_0 = 0$ and $\sigma_0 = \Vert g \Vert$), for $f$
\begin{equation}
 \label{minf}
f_{k+1}
=  {\rm arg min}_{f }\left[ \frac{\| Af-g \|^2}{\sigma_k} + \sigma_k + 2\mu\, \Vert f\Vert_1\right] ,
\end{equation}
or equivalently
\begin{equation}
 \label{minf2}
f_{k+1}
=  {\rm arg min}_{f }\left[ \| Af-g \|^2 + 2\mu \sigma_k\, \Vert f\Vert_1\right], \quad k=0,1, \dots  ,
\end{equation}
and for the residual
\begin{equation}
\label{minsig}
\sigma_{k+1}
=  {\rm arg min}_{\sigma>0}\left[ \frac{\| Af_{k+1}-g \|^2}{\sigma} + \sigma + 2\mu\, \Vert f_{k+1}\Vert_1\right]  ,
\end{equation}
which simply yields
\begin{equation}
 \label{minsig2}
\sigma_{k+1}
=  \| Af_{k+1}-g \|  .
\end{equation}

Now, the idea is to replace the minimisation (\ref{minf}), namely of $\Phi_{\mu}(f, \sigma_k)$ with respect to $f$, which is a Lasso problem, by the minimisation of the surrogate cost function 
\begin{equation}
\label{surk}
{\Phi}^{\scriptstyle\rm{SUR}}_{\mu}(f, \sigma_k; f_k) = \Phi_{\mu}(f, \sigma_k) + \frac{1}{\sigma_k}\left[\frac{1}{\tau} \; \|f-f_k\|^2-\|Af-Af_k\|^2\right]
\end{equation}
where $\tau$ is a positive parameter.  
If $ \|A^*A \| = \Vert A\Vert^2 \leq 1/\tau$, where $\Vert A\Vert$ denotes the operator spectral norm and $A^*$ the adjoint of $A$, this surrogate clearly satisfies the majorisation property
\begin{equation}
\label{major}
\Phi_{\mu}(f,\sigma_k) \le {\Phi}^{\scriptstyle\rm{SUR}}_{\mu}(f, \sigma_k; f_k) , ~~~~  \mbox{{\rm for all}}~f, \sigma_k ~\mbox{{\rm and}} ~f_k .
\end{equation}
It also satisfies the anchoring property
\begin{equation}
\label{anchor}
\Phi_{\mu}(f_k,\sigma_k) = {\Phi}^{\scriptstyle\rm{SUR}}_{\mu}(f_k, \sigma_k; f_k),  ~~  \mbox{{\rm for all}}~f_k~\mbox{{\rm and}} ~\sigma_k  .
\end{equation}
A surrogate like (\ref{surk}) has been used in \cite{Daubechies2004} to derive the ISTA algorithm for the Lasso problem. Optimisation strategies of this type
are usually referred to as MM optimisation (for Majorisation-Minimisation) (see e.g. \cite{Lange2016}).
Thanks to the decoupling of the components of $f$, the minimiser of (\ref{surk}) is easily calculated by hand and is given by
\begin{equation} \label{iterates}
f_{k+1}={S}_{2\tau\mu\sigma_k}\left( f_{k} + \tau A^* (g-Af_{k})\right),
\end{equation}
where ${S}_t$ denotes  the (nonlinear) ``soft-thresholding'' operator, 
acting component-wise, meaning that its action on $h = \sum_\gamma h_\gamma\varphi_\gamma$ is defined by
${S}_t(h) = \sum_\gamma {S}_t(h_\gamma) \varphi_\gamma$, and
where  $S_t(x)$ is the 
soft-thresholding function defined as 
\begin{equation} \label{soft_th}
S_{t}(x)= \left\{
\begin{array}{ccl} x +t/2 ~&~ \mbox{if} ~& x \leq - t/2 \\
0 ~&~ \mbox{if} ~& |x| < t/2
\\ x- t /2 ~&~ \mbox{if} ~& x \geq t/2.
\end{array} \right.  
\end{equation}

The update on $\sigma$ will be as above given by (\ref{minsig2}). The positive parameter $\tau$,  also referred to as a relaxation parameter, characterises the stepsize taken in the direction of the negative gradient of $\Vert Af-g\Vert^2$. 

This iterative algorithm is a simple modification of ISTA, with a threshold varying at each iteration according to the current value $\sigma_k$ of the residual. 
Indeed, for the Lasso cost function (\ref{phi_Lasso}), if the soft-thresholding operator ${S}_{2\tau\mu\sigma_k}$ in (\ref{iterates}) is replaced by ${S}_{\tau\tilde\mu}$ we get the ISTA algorithm
\begin{equation}
f_{k+1}={S}_{\tau\tilde\mu}\!\left( f_{k} + \tau A^* (g-Af_{k})\right),
\end{equation}
(starting from any $f_0$), which can be shown to converge to a minimiser of (\ref{phi_Lasso}) \cite{Daubechies2004}.

Let us remark that a MM-optimisation strategy does not necessarily guarantee the convergence of the sequence of iterates to a minimum point of the cost function. This requires additional work, specific for each case.  

Hence, to summarise, in the next subsection, we will investigate the convergence of the iterative algorithm with iterates $(f_k, \sigma_k)$ and the update steps (\ref{iterates}) and (\ref{minsig2}), namely
\begin{equation} \label{iterates_tau_summ}
f_{k+1}={S}_{2\tau\mu\sigma_k}\!\left( f_{k} + \tau A^* (g-Af_{k})\right), \qquad k=0,1, \dots,
\end{equation}\
\begin{equation}
\label{minsig2_summ}
\sigma_{k+1}
=  \| Af_{k+1}-g \| ,
\end{equation}
initialised with arbitrary $f_0$ such that $\sigma_0 = \Vert Af_0-g\Vert  > 0$. We will refer to this algorithm as SQRT-ISTA.

\begin{remark}
The SQRT-ISTA algorithm (\ref{iterates_tau_summ}, \ref{minsig2_summ}) can also be interpreted as a variable step-length
proximal-gradient algorithm (Theorem 3.4 in \cite{Combettes2005}) applied
to the minimisation of (\ref{phi}). As the gradient of the first term in (\ref{phi}) is
$A^*(Af-g)/\|Af-g\|$, the proximal gradient algorithm for (\ref{phi}) takes
the form
\begin{equation}
f_{k+1}= S_{2\mu \gamma_k}\left(f_k -\gamma_k
\frac{A^*(Af_k-g)}{\|Af_k-g\|}\right).
\end{equation}
Setting $\gamma_k=\tau \|Af_k-g\|$ one formally obtains algorithm   (\ref{iterates_tau_summ}, \ref{minsig2_summ}). However,
the usual convergence theorem (Theorem 3.4 in \cite{Combettes2005}) cannot
be invoked because the first term in (\ref{phi}) is not everywhere
differentiable (not when $Af=g$), its gradient is not Lipschitz
continuous and the stepsizes $\gamma_k$ depend on the iterates
themselves. 

On the other hand, the proximal gradient algorithm of \cite{Fest2023} is related to problem
(\ref{minf2}) in the sense that it updates the penalty parameter in
every step (as in the cost function in expression (\ref{minf2})). However, due to
the dependence of $\sigma_k$ on $f_k$, the algorithm of \cite{Fest2023}
cannot either be applied to the problem discussed here.

Therefore, due to these features, a careful and
separate convergence analysis of the algorithm is necessary.
\end{remark}

\subsection{Convergence properties} 
\label{convergence}

A preliminary observation is in order. Looking at (\ref{iterates_tau_summ}, \ref{minsig2_summ}), we see that when,  for some iteration index $k$, we have $\sigma_k = 0$, 
then $f_k= f_{k+1} = f_{k+2}= \dots$ and the algorithm trivially converges (but not necessarily to a minimiser of (\ref{phi}), as we will see later). Therefore, we can focus only on the case where $\sigma_k >0$ for all $k$, which we assume in the following. 
Notice that the vanishing of $\sigma_k$ for some finite $k$ may result from a mere numerical accident. In such case, one can always try to restart the algorithm with another initialisation. 

A first nice property, typical of  MM-optimisation methods, results from the derivation of the algorithm through surrogates. Indeed, this ensures a monotonic decrease of the cost function $\Phi_{\mu}$ at each iteration and, therefore, ensures its convergence.

\begin{proposition} 
\label{cvcost}
If $\tau \leq 1/\Vert A \Vert^2$, we have along the iteration (\ref{iterates_tau_summ}, \ref{minsig2_summ})
\begin{equation}
\Phi_{\mu}(f_{k+1}) \leq \Phi_{\mu}(f_k); \qquad k=0,1,\dots.
\end{equation} Hence the nonnegative and nonincreasing sequence $\{\Phi_{\mu}(f_k)\}$ has a limit.
\end{proposition}
{\em Proof:}
We have
\begin{equation}
\label{monot1}
2\Phi_\mu(f_{k+1})=\Phi_{\mu}(f_{k+1}, \sigma_{k+1}) \leq \Phi_{\mu}(f_{k+1}, \sigma_{k}) 
\end{equation}
since $\sigma_{k+1}$ minimises $\Phi_{\mu}(f_{k+1}, \sigma)$. Moreover
\begin{equation}
 \label{monot2}
\begin{array}{lcl}
\Phi_{\mu}(f_{k+1}, \sigma_{k}) &\leq&{\Phi}^{\scriptstyle\rm{SUR}}_{\mu}(f_{k+1}, \sigma_k; f_k)\\[3mm] 
&\leq& {\Phi}^{\scriptstyle\rm{SUR}}_{\mu}(f_{k}, \sigma_k; f_k) = \Phi_{\mu}(f_{k}, \sigma_{k}) =2 \Phi_\mu(f_k)
\end{array}
\end{equation}
by the fact that $f_{k+1}$ minimises (\ref{surk}) and by the majorisation and anchoring properties (\ref{major}) and (\ref{anchor}) of the surrogate.
\hfill \QED

We will now establish convergence properties of the iterates in a slightly generalised setting.  We no longer require that $ \|A^*A \| < 1/\tau$, hence losing the majorisation property (\ref{major}) but, as we will see, this will allow to prove convergence under a more flexible assumption on the link between the stepsize $\tau$ along the negative gradient and the norm of $A$. 

To analyse the convergence of the iterative scheme SQRT-ISTA (\ref{iterates_tau_summ}, \ref{minsig2_summ}), we first establish two lemmas.
The first lemma will enable us to control the decrease of the cost function for two successive iterates.

\begin{lemma} \label{descent_disc}
The following inequality holds true for any pair of vectors $f_{k+1}$ and $f_{k}$
\begin{equation} \label{ineq_discrep}
\| Af_{k+1}-g \|^2 - \| Af_{k} -g \|^2 - 2(A^* (Af_{k}-g), f_{k+1} - f_{k}) \le \| A \|^2 \| f_{k+1} -f_{k}  \|^2.
\end{equation} 
\end{lemma} 
{\em Proof:}  The inequality follows by noting that, using elementary manipulations,  the left-hand side of (\ref{ineq_discrep})
is equal to $\|A( f_{k+1} -f_{k})  \|^2$.
 \hfill \QED

Let us recall the definition of the ``prox'' $\bar f$ of $h$ with respect to a function $\Psi(f)$ (here the penalty)
\begin{equation}
\bar f = {\rm arg min}_f \left[ \frac{1}{2}\  \Vert f-h\Vert^2 + \Psi (f) \right].
\end{equation} 
It is a well-known fact that for $\Psi (f)= \mu\sigma \Vert f \Vert_1$ 
we have $\bar f = S_{2\mu\sigma} (h)$ \cite{Beck2017}.

We will repeatedly use the subdifferential of a convex function $\Psi:\mathcal{F}\to \mathbb{R}$ in a point $f$, which is defined as
\begin{equation}\label{eq:subdiff}
\partial \Psi (f) = \{w \in \mathcal{F}\ {\rm such\ that }\ \Psi(\tilde f)\geq \Psi(f)+(w,\tilde f-f)\ {\rm for\ all}\ \tilde f \in \mathcal{F}\}.
\end{equation}
Its principal property is that $0\in \partial \Psi(f)$ if and only $f$ is a minimiser of $\Psi$. The elements of the subdifferential are called subgradients.
With this definition, the second lemma is standard in convex analysis and subdifferential calculus (see \cite{Bauschke2011}, Section 16.4).
\begin{lemma} \label{optsur2}
Let $h\in \mathcal{F}$. The relation $\bar f= S_{2\tau\mu\sigma} (h) = {\rm arg min}_f  \left( \Vert f-h\Vert^2 + 2\tau\mu\sigma \Vert f \Vert_1\right)$ is equivalent to the inclusion  
\begin{equation}\label{subdiff-soft}
h-\bar f \in  \partial(\tau\mu\sigma \Vert \bar f \Vert_1 )
\end{equation}
and to the inequalities (for all $f\in \mathcal{F}$)
\begin{equation} \label{subdiff2}
\tau \mu\sigma \Vert \bar f \Vert_1 \leq \tau \mu\sigma \Vert f \Vert_1 -   (h-\bar f, f- \bar f).
\end{equation}
\end{lemma} 
{\em Proof:} $\bar f$ is the mimimizer of $\Vert f-h\Vert^2 + 2\tau\mu\sigma \Vert f \Vert_1$ if and only if $0$ is an element of its subdifferential evaluated at $\bar f$, i.e. $0 \in 2(\bar f-h) +\partial(2\tau \mu\sigma \Vert \bar f \Vert_1 )$ or 
\begin{equation} \label{subdiffeq}
2(h-\bar f) \in  \partial(2\tau \mu\sigma \Vert \bar f \Vert_1 ).
\end{equation}
By definition of the subdifferential (\ref{eq:subdiff}), this corresponds to the inequalities (\ref{subdiff2}). 
\hfill \QED

Combining these two lemmas, we can establish the monotonic decrease of the cost function along the iteration, generalising Proposition \ref{cvcost}.

\begin{proposition} \label{decay_Phi} Provided that $0 < \tau \leq 2/ \| A \|^2$,
the sequence of values of the cost function is monotonically decreasing along the iteration (\ref{iterates_tau_summ}, \ref{minsig2_summ}):
\begin{equation} \label{decrease_cost}
\Phi_{\mu}(f_{k+1}) \leq \Phi_{\mu}(f_k), \qquad k=0, 1, \dots.
\end{equation}
Hence the nonnegative and nonincreasing sequence $\{\Phi_{\mu}(f_k)\}$ has a limit.
When $0<\tau<2/\|A\|^2$, the decrease is strict, i.e. $\Phi_\mu (f_{k+1}) < \Phi_\mu (f_{k})$ unless $f_{k+1} = f_{k}$.
\end{proposition}
{\em Proof:}
We consider the inequality (\ref{subdiff2}) of Lemma \ref{optsur2} for $f=f_k$, $h= f_{k} + \tau A^* (g-Af_{k})$ and, according to (\ref{iterates_tau_summ}), $\bar f=f_{k+1}$, namely
\begin{equation} \label{decay_penalty}
\tau\mu\sigma_k \Vert f_{k+1} \Vert_1 \leq \tau\mu\sigma_k \Vert f_{k} \Vert_1 -   (f_{k} + \tau A^* (g-Af_{k}) - f_{k+1}, f_{k} - f_{k+1}).\end{equation}
Adding to this latter inequality -- multiplied by $2$ and divided by $\tau$ -- the inequality (\ref{ineq_discrep}) of Lemma \ref{descent_disc}, then dividing by $\sigma_{k} > 0$ and adding $\sigma_{k}$ to both sides, we get, using (\ref{phisig}), 
\begin{equation}
\Phi_\mu(f_{k+1}, \sigma_{k}) \leq \Phi_\mu(f_{k}, \sigma_{k}) + \frac{1}{\sigma_{k}} \left(\| A \|^2 - \frac{2}{\tau}\right) \| f_{k+1} - f_{k} \|^2.
\end{equation}
Using the fact that $\sigma_{k+1}$ minimises $\Phi_\mu(f_{k+1}, \sigma)$, we finally get
\begin{equation} \label{decay_bound}
0 \leq \Phi_\mu(f_{k+1}, \sigma_{k+1}) \leq \Phi_\mu(f_{k}, \sigma_{k}) + \frac{1}{\sigma_{k}} \left(\| A \|^2 - \frac{2}{\tau}\right) \| f_{k+1} - f_{k} \|^2
\end{equation} 
or
\begin{equation} \label{decay_bound2}
0 \leq \Phi_\mu(f_{k+1}) \leq \Phi_\mu(f_{k}) + \frac{1}{2\sigma_{k}} \left(\| A \|^2 - \frac{2}{\tau}\right) \| f_{k+1} - f_{k} \|^2.
\end{equation} 
 Hence, provided the last term is not positive, we can deduce that
\begin{equation} \label{decrease_cost2}
\Phi_\mu(f_{k+1}) \leq \Phi_\mu(f_{k}) \leq \dots \leq \Phi_\mu(f_{0}).
\end{equation}
When $0<\tau<2/\|A\|^2$, we also see from (\ref{decay_bound2}) that the decrease of cost function is strict, i.e. $\Phi_\mu (f_{k+1}) < \Phi_\mu (f_{k})$ unless $f_{k+1} = f_{k}$.  \hfill \QED

This in turn implies that the sequences of the iterates and of the residuals are bounded.
\begin{proposition} \label{prop_bddness} Provided that $0< \tau \leq 2/ \| A \|^2$,
the sequence of iterates $\{(f_k,\sigma_k)\}$ is bounded.
\end{proposition}
{\em Proof:}
We have by (\ref{phisig}) and (\ref{decrease_cost})
\begin{equation} \label{phibdd}
\Phi_\mu (f_{k}) = \sigma_{k} + \mu \Vert f_{k} \Vert_1 \leq \Phi_\mu(f_{0}).
\end{equation}
Since $\mu$ is fixed, we see that $ \Vert f_k \Vert_1$ is uniformly bounded in $k$. By the inequality
\begin{equation}
\Vert f \Vert^2 = \sum_\gamma f_\gamma^2 \leq  {\rm max}_\gamma \{\vert f_\gamma\vert\}
\sum_\gamma  \vert f_\gamma\vert \leq  (\sum_\gamma \vert f_\gamma\vert )^2 = \Vert f \Vert_1^2,
\end{equation}
we see that the norms $\Vert f_k\Vert$ are bounded as well. 

Notice that this remains true in the case of a weighted $\ell^1$-norm, since with $c$ the strictly positive lower bound of the weights, we have $\Vert f \Vert^2 = \sum_\gamma f_\gamma^2 \leq  {\rm max}_\gamma \left\{ \frac{w_\gamma}{w_\gamma^2}\ \vert f_\gamma\vert\right\} \sum_\gamma w_\gamma \vert f_\gamma\vert \leq \frac{1}{c^2} ( \sum_\gamma w_\gamma \vert f_\gamma\vert)^2 = \frac{1}{c^2} \Vert f \Vert_1^2$. 

Moreover, by (\ref{phibdd}), we see that the residuals $\sigma_k$ are also uniformly bounded.  
\hfill \QED

Using these results, we can now establish the so-called property of asymptotic regularity of the iterative scheme. This means that the norm of the difference between two successive iterates tends to zero with the number of iterations.

\begin{proposition} \label{asreg} (``asymptotic regularity'')  
Provided that $0< \tau <  2/ \| A \|^2$, for two successive iterates $f_k$ and $f_{k+1}$, we have that $\Vert f_{k+1} - f_k \Vert \to 0$ for $k \to \infty$. For the corresponding residuals, we also have $\vert \sigma_{k+1} - \sigma_k \vert \to 0$.
\end{proposition}
{\em Proof:}
Let us denote by $\sigma_{\max}$ an upper bound for all $\sigma_{k}$'s. Since $\sigma_{k} \leq \sigma_{\max}$, we have by (\ref{decay_bound})
\begin{eqnarray} \label{decay_cost_asreg}
\frac{1}{\sigma_{\max}} \left( \frac{2}{\tau} - \| A \|^2 \right) \| f_{k+1} - f_{k} \|^2 &\leq& \frac{1}{\sigma_{k}} \left( \frac{2}{\tau} - \| A \|^2 \right)  \| f_{k+1} - f_{k} \|^2 \nonumber \\ &\leq&  \Phi_\mu (f_{k}, \sigma_{k})
-\Phi_\mu (f_{k+1}, \sigma_{k+1}). \end{eqnarray}
Hence, by a telescopic sum argument,
\begin{eqnarray}
\frac{1}{\sigma_{\max}} \left( \frac{2}{\tau} - \| A \|^2 \right) \sum_{k=0}^K \| f_{k+1} - f_{k} \|^2 &\leq& \sum_{k=0}^K \left[\Phi_\mu(f_{k}, \sigma_{k}) - \Phi_\mu(f_{k+1}, \sigma_{k+1})\right] \nonumber\\
&=&  \Phi_\mu(f_{0}, \sigma_{0}) - \Phi_\mu(f_{K+1}, \sigma_{K+1}) \nonumber\\
 &\leq& \Phi_\mu(f_{0}, \sigma_{0}).
\end{eqnarray}
From this we conclude that
\begin{equation}
\sum_{k=0}^\infty \| f_{k+1} - f_{k} \|^2 < + \infty, \ \ \hbox{implying} \ \  \| f_{k+1} - f_{k} \| \to 0\ \  \hbox{for}\ \  k \to \infty.
\end{equation}
Moreover,
\begin{eqnarray} \label{asregsigma}
 \vert \sigma_{k+1} - \sigma_k \vert &=& \vert \  \Vert Af_{k+1} - g\Vert -  \Vert Af_{k} - g\Vert\  \vert \nonumber \\
&\leq& \Vert A(f_{k+1} - f_k) \Vert \leq \Vert A \Vert \  \Vert f _{k+1} - f_k \Vert
\end{eqnarray}
using the reverse triangular inequality. Hence the asymptotic regularity for the residuals follows from that of the iterates $f_k$.
\hfill \QED

To proceed further, we exploit the fact that the sequence of iterates is bounded. Therefore, in finite dimension, it contains convergent subsequences. As we will see, this last property is crucial for our convergence proof. Accordingly, although the derivation up to this point would hold also in an infinite-dimensional Hilbert space setting, we will assume from now on that $f$ is a vector belonging to $\mathbb R^d$. We will also assume that $0 <\tau < 2/\Vert A \Vert^2$.

Let us denote by $T$ the combined iteration mapping that maps the pair of iterates $(f_k, \sigma_k)$ onto the next one $(f_{k+1}, \sigma_{k+1})$. According to (\ref{iterates_tau_summ}, \ref{minsig2_summ}), it is given by 
\begin{eqnarray}\label{itermap}
T(f, \sigma) = (\tilde f, \tilde \sigma) \ \hbox{with}\ &\tilde f&={S}_{2\tau\mu\sigma}\!\left( f + \tau A^* (g-Af)\right),\nonumber\\ 
 &\sigma& = \Vert A f - g\Vert \ \hbox{and}\ \tilde \sigma = \Vert A \tilde f - g\Vert. 
\end{eqnarray}
We observe that $T$ is a continuous mapping in both $f$ and $\sigma$. 
Notice that contrary to the mapping applying $f_k$ on $f_{k+1}$, which depends on $k$ via the threshold $\sigma_k$, this combined mapping is stationary, i.e. does not depend on $k$. The next proposition shows that the limit $(f_*, \sigma_*)$ of any convergent subsequence of the (bounded) sequence of iterates is a fixed point of the mapping $T$, namely that $(f_*,\sigma_*)= T(f_*,\sigma_*)$.

\begin{proposition}\label{fixedpt}
Let $\{(f_{k_j},\sigma_{k_j})\}$ be a convergent subsequence of the sequence of iterates $\{(f_k,\sigma_k)\}$, with limit $(f_*, \sigma_*)$. Then it converges to a fixed point of the iteration mapping $T$, i.e. $(f_*,\sigma_*)= T(f_*,\sigma_*)$, or equivalently
\begin{equation} \label{fixedpteq}
f_*= {S}_{2\tau\mu\sigma_*}\!\left( f_* + \tau A^* (g-Af_*)\right) \ \hbox{with} \  \sigma_*= \Vert A f_* - g\Vert.
\end{equation}
\end{proposition}
{\em Proof:}  If $(f_{k_j},\sigma_{k_j})\stackrel{j\to
\infty}{\longrightarrow}(f_*,\sigma_*)$ then, by
Proposition \ref{asreg}, also $(f_{k_j+1},\sigma_{k_j+1})\stackrel{j\to
\infty}{\longrightarrow}(f_*,\sigma_*)$. By taking the same limit
on both sides of the update rule
$(f_{k_j+1},\sigma_{k_j+1})=T(f_{k_j},\sigma_{k_j})$, one finds, using
the continuity of $T$, that $(f_*,\sigma_*)=T(f_*,\sigma_*)$.
\hfill \QED

\begin{proposition}
\label{minpt}
 Provided that $\sigma_* \neq 0$, a fixed point of the iteration $(f_*,\sigma_*)= T(f_*,\sigma_*)$ with $T$ given by (\ref{itermap}) is a minimiser of the cost function (\ref{phi}).
\end{proposition}
{\em Proof:} When $\sigma_\ast=\|Af_\ast-g\|>0$, the gradient of $\|Af-g\|$ at $f=f_\ast$ is $A^\ast(Af_\ast -g)/\sigma_\ast$ and the optimality conditions for minimising the cost function (\ref{phi}) are
\begin{equation}
0\in A^\ast(Af_\ast -g)/\sigma_\ast +\partial \left(\mu \|f_\ast\|_1\right)
\end{equation}
or 
\begin{equation}
\tau A^\ast(g-Af_\ast )\in  \partial \left(\tau \mu \sigma_\ast\|f_\ast\|_1\right)
\end{equation}
for any $\tau>0$.
By Lemma~\ref{optsur2} and equation (\ref{subdiff-soft}) with $h=f_\ast+\tau A^\ast(g-Af_\ast )$ and $\bar f=f_\ast$, this inclusion is equivalent to the relation $f_\ast={S}_{2\tau\mu\sigma_*}\!\left( f_* + \tau A^* (g-Af_*)\right)$.
\hfill \QED

\begin{remark}
Using the optimality conditions, we can see that the minimisers of the Lasso cost function $\|Af-g\|^2+\tilde\mu \Vert f \Vert_1$ and the
minimisers of the square-root Lasso cost function $ \|Af-g\|+\mu
\Vert f \Vert_1$ coincide when $\tilde\mu= 2 \mu \|Af_*-g\|$ since they
respectively satisfy
\begin{equation}
0 \in 2A^*(Af_*-g)+ \tilde\mu \partial \Vert f_* \Vert_1\ \hbox{and}\
0 \in A^*(Af_*-g)/\|Af_*-g\|+ \mu \partial \Vert f_* \Vert_1
\end{equation}
when $Af_*\neq g$.  This was already observed in \cite{Giraud2015} and \cite[Section 3.7]{VandeGeer2016}.  Notice, however, that if the square-root Lasso has multiple solutions for a given $\mu$, then these may correspond to Lasso solutions for different values of $\tilde\mu$.
\end{remark}

Provided that the residuals $\sigma_k$ do not vanish in the limit, Proposition \ref{prop_bddness} implies (in finite dimension) that a converging subsequence must exist, whose limit is a minimiser of (\ref{phi}) (Propositions \ref{fixedpt} and \ref{minpt}). If such a minimiser is unique, then the whole sequence converges to this minimiser.
 
Nevertheless, this uniqueness assumption as well as the exclusion of vanishing residuals may appear as restrictive. Furthermore such assumptions are not required to prove the convergence of the ISTA algorithm for Lasso. We will now show that both restrictions can be lifted when resorting to the use of a powerful general tool in optimisation, namely the use of the so-called Kurdyka-Łojasiewicz (KL) inequality.

Relying on earlier work by Łojasiewicz \cite{Lojasiewicz1963} and Kurdyka \cite{Kurdyka1998}, the use of this inequality was later advocated mainly to address the case of nonconvex \cite{Attouch2010} and nonsmooth \cite{Bolte2007} problems where the cost function to be minimised satisfies the KL property (defined here below).
It appears that the class of such functions is very large and contains lots of practical instances. In particular, it is shown in \cite[Appendix 5]{Bolte2014}, that any $p$-norm, with $p>0$ rational, is semi-algebraic and therefore satisfies the KL property. Hence, inspired by the convergence proof of the Theorem 1 in \cite{Bolte2014} in the context of another algorithm (i.e. PALM - Proximal Alternating Linearized Minimization),  we will exploit the KL property of the cost function (\ref{phi}). This will allow us to prove that the sequence of iterates is a Cauchy sequence and therefore always converges. The setting also requires that the set of the limit points of the sequence of iterates is a compact connected set.  But this follows from the following result about the limit points of an iteration with asymptotic regularity,  known as Ostrowski's theorem (see the book by Ostrowski  \cite{Ostrowski1973} or Proposition 8.2.1 in the book by Lange \cite{Lange1999}).

\begin{proposition}
\label{ostrowski}
Let $\{f_k\}_k$ be a bounded sequence in $\mathbb{R}^d$. If, for some norm $\Vert f \Vert$,
\begin{equation} 
\label{asregseq}
 \lim_{k \rightarrow +\infty} \Vert f_{k+1}-f_{k}\Vert =0,
\end{equation}
then the set of the limit points of $\{f_k\}_k$ is a compact connected set. 
Furthermore, if this set is finite, then it reduces to a single point and the sequence converges.
\end{proposition}

We also need the following Lemma.
\begin{lemma}
\label{ineq_soft}
The following inequality holds for the soft-thresholding operator $S_{t}$
\begin{equation} \label{ineq-soft2}
\|S_{t}(x)-x\|\leq \sqrt{d}\ \frac{t}{2} \quad \mbox{for any}\  x\in\mathbb{R}^d.
\end{equation}
\end{lemma}
{\em Proof:}  This is easy to see in the special case $d=1$. Indeed, according to (\ref{soft_th}), we have
\begin{equation}
\left| S_{t}(x)-x\right| =
\left\{
\begin{array}{lcl}
|0-x| = |x| &\mbox{if}& |x| < t/2\\
|x\pm (t/2)  -x| = t/2& \mbox{if} & |x|\geq  t/2
\end{array}
\right\}
\leq  \frac{t}{2} \quad 
\end{equation}

independently of $x$.
The case $x \in \mathbb{R}^d$ with  $S_{t}$ acting component-wise then goes as follows 
$\left\| S_{t}(x)-x\right\|^2=\sum_{i=1}^d |S_{t}(x_i)-x_i |^2\leq \sum_{i=1}^d (t/2)^2=d\; (t/2)^2$.
\hfill \QED

\begin{proposition}\label{lemma:sigmas} The sequence $\{\sigma_{k+1}/\sigma_k\}$ is bounded.
\end{proposition}
{\em Proof:} 
Let us compute the residual $ \sigma_{k+1}$ using (\ref{iterates_tau_summ}) 
\begin{eqnarray}
\displaystyle \sigma_{k+1} &= & \displaystyle \left\|Af_{k+1}-g\right\| 
 = \displaystyle \left\|A\; S_{2\tau\mu\sigma_k}\!\left( f_{k} + \tau A^* (g-Af_{k})\right)-g\right\| \nonumber \\
&\leq & \displaystyle \left\|A\left( f_{k} + \tau A^* (g-Af_{k})\right)-g\right\| \nonumber \\
&& \displaystyle +\left\|A\; S_{2\tau\mu\sigma_k}\!\left( f_{k} + \tau A^* (g-Af_{k})\right)-A\left( f_{k} + \tau A^* (g-Af_{k})\right)\right\|  \nonumber\\
&\leq & \displaystyle \left\|A  f_{k}-g\right\|+\left\|\tau AA^* (g-Af_{k})\right\|  \nonumber\\
&& \displaystyle +\|A\|\,\left\|S_{2\tau\mu\sigma_k}\!\left( f_{k} + \tau A^* (g-Af_{k})\right)-\left( f_{k} + \tau A^* (g-Af_{k})\right)\right\|  \nonumber\\
&\leq & \displaystyle \sigma_k+\tau \|A\|^2 \left\|g-Af_{k}\right\|+\|A\|\,\sqrt{d}\ \tau\mu\ \sigma_k  \nonumber\\
&\leq & \displaystyle C\ \sigma_k
\end{eqnarray}
with $C \geq 1 $ a constant independent of $k$, where we have also used $\sigma_k=\|Af_k-g\|$ and Lemma \ref{ineq_soft}.
\hfill \QED
We now characterise the distance between the origin and the subdifferential of $ \Phi_\mu$ at the point $f_{k+1}$.
\begin{proposition} We have that 
\begin{equation}\label{subdiffbound}
\dist(0,\partial \Phi_\mu(f_{k+1})) \leq C_1 \|f_{k+1}-f_k\|/\sigma_k
\end{equation}
with $C_1$ a positive constant independent of $k$.
\end{proposition}
{\em Proof:} Rewriting iteration (\ref{iterates_tau_summ}) as $f_{k+1}={S}_{2\tau\mu\sigma_k}\!\left(f_{k+1}+ f_{k}-f_{k+1} - \tau A^* (Af_{k}-g)\right)$, we have by Lemma \ref{optsur2} and equation (\ref{subdiff-soft}) that 
\begin{equation}
f_{k}-f_{k+1} - \tau A^* (Af_{k}-g) \in  \partial\left( \tau\mu\sigma_k  \|f_{k+1}\|_1\right)
\end{equation}
from which it follows
\begin{eqnarray}
\frac{f_{k}-f_{k+1}}{\tau\sigma_k}&-&\frac{A^* (Af_{k}-g)}{\sigma_k} + \frac{A^* (Af_{k+1}-g)}{\sigma_{k+1}}\nonumber\\
&& \in \frac{A^* (Af_{k+1}-g)}{\sigma_{k+1}}+ \partial\left(\mu  \|f_{k+1}\|_1\right).
\end{eqnarray}
As $A^*(Af_{k+1}-g)/\sigma_{k+1}$ is the gradient of $\|Af-g\|$ evaluated at $f=f_{k+1}$, we obtain the inclusion
\begin{equation}
\frac{f_{k}-f_{k+1}}{\tau\sigma_k}-\frac{A^* (Af_{k}-g)}{\sigma_k} + \frac{A^* (Af_{k+1}-g)}{\sigma_{k+1}} \in \partial\Phi_\mu(f_{k+1}).
\end{equation}
Since $\dist(0,\partial \Phi_\mu(f_{k+1}) )$ is bounded by the norm of any element of $\partial \Phi_\mu(f_{k+1})$, we write 
\begin{eqnarray}
&\dist&(0,\partial \Phi(f_{k+1}) )  \leq  \displaystyle \left\| \frac{f_{k}-f_{k+1}}{\tau\sigma_k}-\frac{A^* (Af_{k}-g)}{\sigma_k} + \frac{A^* (Af_{k+1}-g)}{\sigma_{k+1}}\right\| \nonumber\\
 &\leq&\displaystyle \frac{\|f_{k}-f_{k+1}\|}{\tau\sigma_k}+\left\|-\frac{A^* (Af_{k}-g)}{\sigma_k}+\frac{A^* (Af_{k+1}-g)}{\sigma_k}\right\| \nonumber\\
 && \displaystyle \qquad\qquad\qquad+\left\|-\frac{A^* (Af_{k+1}-g)}{\sigma_k} + \frac{A^* (Af_{k+1}-g)}{\sigma_{k+1}} \right\| \nonumber\\
&\leq&\displaystyle \frac{\|f_{k}-f_{k+1}\|}{\tau\sigma_k}+\frac{\|A\|^2 \left\|f_{k}-f_{k+1}\right\|}{\sigma_k}
 \displaystyle +\|A\| \|Af_{k+1}-g\| \left|\frac{1}{\sigma_{k+1}}-\frac{1}{\sigma_k}  \right| \nonumber\\
& \stackrel{(\ref{minsig2_summ})}{=}& \displaystyle \frac{\|f_{k}-f_{k+1}\|}{\tau\sigma_k}+\frac{\|A\|^2 \left\|f_{k}-f_{k+1}\right\|}{\sigma_k} +\|A\| \frac{\left|\sigma_{k+1}-\sigma_k\right|}{\sigma_k}\nonumber\\
& \stackrel{(\ref{asregsigma})}{=}& \displaystyle \frac{\|f_{k}-f_{k+1}\|}{\tau\sigma_k}+\frac{\|A\|^2 \left\|f_{k}-f_{k+1}\right\|}{\sigma_k} +\|A\|^2 \frac{\left\|f_{k+1}-f_k\right\|}{\sigma_k}
\end{eqnarray}
which implies inequality (\ref{subdiffbound}).\hfill \QED

We now have all ingredients to proceed to the proof of convergence of the algorithm. 
 
\begin{theorem} \label{final_cv}
When $0<\tau<2/\|A\|^2$, the sequence $\{(f_k,\sigma_k)\}_k$ defined by the iteration (\ref{iterates_tau_summ},\ref{minsig2_summ}) converges. If $\lim_k\sigma_k>0$, then $\lim_k f_k$ is a minimizer of $\Phi_\mu$.
\end{theorem}
{\em Proof:} We have already established in Proposition \ref{prop_bddness} the boundedness of $\{f_k\}_k$ and $\{\sigma_k\}_k$. We set $f_\ast$ to be the limit of a converging subsequence of $\{f_k\}_k$. We also know by Proposition \ref{decay_Phi} that $\Phi_\mu(f_k)$ decreases (strictly when $0<\tau<2/\|A\|^2$, otherwise we would have $f_{k+1}=f_k$ which trivially converges). Hence $\Phi_\mu(f_k)>\Phi_\mu(f_\ast)$ and $\Phi_\mu(f_k)\stackrel{k\to \infty}{\rightarrow}\Phi_\mu(f_\ast)$ (and not just for a subsequence).

We set $\omega$ to be the set of limit points of the sequence  $\{f_k\}_k$. According to Proposition \ref{ostrowski}, this is  a compact and connected set. 
We then have that $\lim_k \dist(\omega,f_k)=0$ and $\Phi_\mu$ is constant on $\omega$, its value being $\Phi_\mu(f_\ast)$.

Relying on \cite[Lemma 6]{Bolte2014}, we now use the fact that the cost function $\Phi_\mu$ has the uniformised Kurdyka-Łojasiewicz property, meaning there exist $\epsilon, \eta>0$ and a function $\varphi: [0, \eta[\to \mathbb{R}_{\geq 0}$
with the following properties:
\begin{itemize}
\item $\varphi$ is concave;
\item $\varphi$ is continuous on $[0,\eta[$;
\item $\varphi$ is differentiable on $]0,\eta[$;
\item $\varphi(0)=0$;
\item $\varphi'>0$ on $]0,\eta[$;
\item the Kurdyka-Łojasiewicz inequality:
\begin{equation}\label{KLineq}
\varphi'(\Phi_\mu(f)-\Phi_\mu(f_\ast)) \cdot \dist(0,\partial \Phi_\mu(f) ) \geq 1
\end{equation}
is valid on $\{f \in \mathbb{R}^d \mbox{ such that } \dist(\omega, f)<\epsilon \mbox{ and } \Phi_\mu(f_\ast)<\Phi_\mu(f)<\Phi_\mu(f_\ast)+\eta\}$.
\end{itemize}
Now, as $\Phi_\mu(f_k)\stackrel{k\to \infty}{\rightarrow}\Phi_\mu(f_\ast)$ and $\lim_k \dist(\omega,f_k)=0$, we can use inequality (\ref{KLineq}) with $f=f_k$ for $k$ larger than some $k_{\min}$. For simplicity of notation and without loss of generality, we will assume that $k_{\min}=0$, i.e. we have that
\begin{equation}
\varphi'(\Phi_\mu(f_k)-\Phi_\mu(f_\ast)) \cdot \dist(0,\partial \Phi_\mu(f_k) ), \geq 1\qquad \forall k\geq 0.
\end{equation}
We now combine this inequality with inequality (\ref{subdiffbound}) to obtain
\begin{equation}
\varphi'(\Phi_\mu(f_k)-\Phi_\mu(f_\ast))\geq \frac{1}{C_1}\; \sigma_{k-1}/\|f_k-f_{k-1}\|,
\end{equation}
and express the fact that the concave function $\varphi$ lies below its tangents:
\begin{eqnarray}
&&\varphi(\Phi_\mu(f_{k+1})-\Phi_\mu(f_\ast)) \leq \displaystyle \varphi(\Phi_\mu(f_{k})-\Phi_\mu(f_\ast)) \nonumber\\
&&\displaystyle \qquad\qquad\qquad\qquad\quad\ \  -\varphi'(\Phi_\mu(f_{k})-\Phi_\mu(f_\ast))\cdot \left(\Phi_\mu(f_{k})-\Phi_\mu(f_{k+1})\right)\nonumber\\
&&\leq \displaystyle\varphi(\Phi_\mu(f_{k})-\Phi_\mu(f_\ast)) - \frac{1}{C_1}\, \frac{\sigma_{k-1}}{\|f_k-f_{k-1}\|}\cdot \left(\Phi_\mu(f_{k})-\Phi_\mu(f_{k+1})\right)\nonumber\\
&&\stackrel{(\ref{decay_bound2})}{\leq} \displaystyle\varphi(\Phi_\mu(f_{k})-\Phi_\mu(f_\ast)) - C_2\, \frac{\sigma_{k-1}}{\|f_k-f_{k-1}\|}\frac{\|f_{k+1}-f_k\|^2}{\sigma_k}
\end{eqnarray}
with $C_2>0$, using $\Phi_\mu(f_{k})>\Phi_\mu(f_{k+1})$. 
Setting $a_k=\varphi(\Phi_\mu(f_{k})-\Phi_\mu(f_\ast))>0$ (and likewise $a_{k+1}$), we find by Proposition~\ref{lemma:sigmas}:
\begin{eqnarray}
\|f_{k+1}-f_k\|^2 &\leq& \frac{1}{C_2}\; \frac{\sigma_k}{\sigma_{k-1}}\|f_k-f_{k-1}\|\, \left( a_{k}- a_{k+1}\right)\nonumber\\ &\leq& C_3\; \|f_k-f_{k-1}\|\, \left( a_{k}- a_{k+1}\right).
\end{eqnarray}
where $C_3= C/C_2$.
Using the inequality $2\sqrt{ab}\leq a+b$ for positive numbers $a$ and $b$, we can now obtain
\begin{eqnarray}
2 \|f_{k+1}-f_k\| =  2 \sqrt{\|f_{k+1}-f_k\|^2} &\leq&  2\sqrt{C_3 \|f_k-f_{k-1}\|\, \left( a_{k}- a_{k+1}\right)} \nonumber\\ &\leq &  \|f_k-f_{k-1}\|+ C_3\left( a_{k}- a_{k+1}\right)
\end{eqnarray}
and therefore 
\begin{eqnarray}
2&\sum_{k=1}^{K}&\|f_{k+1}-f_k\| \leq \sum_{k=1}^{K}\left[ \|f_k-f_{k-1}\|+ C_3\left( a_{k}- a_{k+1}\right)\right]\nonumber\\
 &=&\|f_0-f_1\|+\sum_{k=1}^{K} \|f_{k+1}-f_{k}\|-\|f_{K+1}-f_K\|+ C_3\left( a_{1}- a_{K+1}\right) \nonumber\\
 &\leq & \|f_0-f_1\|+\sum_{k=1}^{K}\|f_{k+1}-f_{k}\|+ C_3\ a_{1}
\end{eqnarray}
so that
\begin{equation}
\sum_{k=1}^{K}\|f_{k+1}-f_k\| \leq \|f_0-f_1\|+ C_3\ a_{1}
\end{equation}
independently of $K$. Therefore $\sum_{k=1}^{\infty}\|f_{k+1}-f_k\|<\infty$. This implies that $\{f_k\}_k$ is a Cauchy sequence:
\begin{equation}
\|f_p-f_q\| \stackrel{q>p}{=}\left\|\sum_{k=p}^{q-1} f_k-f_{k+1}\right\|\leq \sum_{k=p}^{q-1} \left\|f_k-f_{k+1}\right\|\stackrel{p,q\to\infty}{\rightarrow}0.
\end{equation}
Hence $\{f_k\}_k$ converges. 

If we set $f_\ast=\lim_{k\to\infty}f_k$ then $\lim_{k\to\infty}\sigma_k=\lim_{k\to\infty}\|Af_k-g\|=\|Af_\ast-g\|$ which we call $\sigma_\ast$. Taking the limit $k\to\infty$ in the relations (\ref{iterates_tau_summ}) and (\ref{minsig2_summ}), one obtains (\ref{fixedpteq}) which, when $\sigma_\ast>0$, by Proposition \ref{minpt}, guarantees that $f_\ast$ is a minimizer of $\Phi_\mu$.\hfill \QED

To assess how fast the iterates converge, we can derive the following convergence rate.
\begin{proposition}\label{proposition:rate} Let $\{(f_k,\sigma_k)\}_k$ be defined by the iteration (\ref{iterates_tau_summ},\ref{minsig2_summ}) and $f_\ast =\lim_{k}f_k$. When $0<\tau\leq 1/\|A\|^2$ and $\sigma_k\geq\sigma_{\min}>0$ for all $k$, the following convergence rate applies:
\begin{equation}
\Phi_\mu(f_k)-\Phi_\mu(f_\ast)\leq \|f_0-f_\ast\|^2/(2\tau\sigma_{\min}k).
\end{equation} 
\end{proposition}
{\em Proof:} By Lemma \ref{optsur2}, iteration (\ref{iterates_tau_summ}) implies $f_k-f_{k+1}-\tau A^\ast (Af_k-g)\in \partial \left(\tau \mu \sigma_k \|f_{k+1}\|_1\right)$ which in turn implies
\begin{equation}
\tau \mu \sigma_k \|f\|_1 \geq \tau \mu \sigma_k \|f_{k+1}\|_1+ \left( f_k-f_{k+1}-\tau A^\ast (Af_k-g), f-f_{k+1} \right)
\end{equation}
for all $f$, or
\begin{eqnarray}
\left( f_k-f_{k+1}, f-f_{k+1} \right) & \leq& \left(\tau A^\ast (Af_k-g), f-f_{k+1} \right) \nonumber\\ &+&\tau \mu \sigma_k \left(\|f\|_1-\|f_{k+1}\|_1\right),
\end{eqnarray}
which we rewrite after simple manipulations as
\begin{eqnarray}
\|f_{k+1}-f\|^2 &\leq&\displaystyle \|f_k-f\|^2 -\|f_{k+1}-f_k\|^2+ 2\tau \left(  A^\ast (Af_k-g), f-f_{k}\right) \nonumber\\
 & +& 2\tau \left(  A^\ast (Af_k-g), f_k-f_{k+1}\right) +2\tau \mu \sigma_k \left(\|f\|_1-\|f_{k+1}\|_1\right).
\end{eqnarray}
The gradient of $\|Af-g\|$ at $f=f_k$ is $A^\ast (Af_k-g)/\sigma_k$. It then follows from the convexity of that function that 
\begin{equation}
\left( A^\ast(Af_k-g)/\sigma_k,f-f_k \right) \leq \|Af-g\|-\|Af_k-g\|
\end{equation}
(for all $f$) yielding:
\begin{eqnarray}
\|f_{k+1}-f\|^2 &\leq&\displaystyle \|f_k-f\|^2 -\|f_{k+1}-f_k\|^2 \nonumber\\
&+& 2\tau\sigma_k \|Af-g\|-2\tau\sigma_k\|Af_k-g\| \nonumber\\
 &+& 2\tau \left(  A^\ast (Af_k-g), f_k-f_{k+1}\right) +2\tau \mu \sigma_k \left(\|f\|_1-\|f_{k+1}\|_1\right).
\end{eqnarray}
Now we use the elementary inequality $\|a\|/2+\|b\|^2/(2\|a\|)\geq \|b\|$, which we rewrite as $\|a\|/2+\|(b-a)+a\|^2/(2\|a\|)\geq \|b\|$ to find:
$\|b\|-\|a\|-( b-a, a/\|a\|)\leq \|b-a\|^2/(2\|a\|)$.
Choosing $a=Af_k-g$ and $b=Af_{k+1}-g$, this gives:
\begin{eqnarray}
&&\|Af_{k+1}-g\|-\|Af_k-g\|-( A(f_{k+1}-f_k), Af_k-g)/\sigma_k \nonumber\\ && \qquad \leq \|A(f_{k+1}-f_k)\|^2/(2\sigma_k)
  \leq  \|A\|^2 \|f_{k+1}-f_k\|^2/(2\sigma_k)
\end{eqnarray}
Inserting this bound in the inequality above yields
\begin{eqnarray}
&& \|f_{k+1}-f\|^2 \leq  \|f_k-f\|^2 -(1-\tau\|A\|^2)\|f_{k+1}-f_k\|^2\nonumber\\
&&\qquad +2\tau\sigma_k \|Af-g\|-2\tau\sigma_k \|Af_{k+1}-g\|+2\tau\mu \sigma_k \left(\|f\|_1-\|f_{k+1}\|_1\right)
\end{eqnarray}
or
\begin{eqnarray}\label{tmp1}
&&\|f_{k+1}-f\|^2 \leq \|f_k-f\|^2 -(1-\tau\|A\|^2)\|f_{k+1}-f_k\|^2\nonumber\\&& \qquad+2\tau \sigma_k \left(\Phi_\mu(f)-\Phi_\mu(f_{k+1})\right).
\end{eqnarray}
Setting $f=f_\ast=\lim_k f_k$ (a minimiser of $\Phi_\mu$) and using $\sigma_k\geq \sigma_{\min}>0$, we find
\begin{equation}\label{eq:descent}
\Phi_\mu(f_{k+1})-\Phi_\mu(f_\ast) \leq \left(\|f_k-f_\ast\|^2-\|f_{k+1}-f_\ast\|^2\right)/(2\tau \sigma_{\min}).
\end{equation}
Finally we use the monotonicity of $\{\Phi_\mu(f_k)\}_k$ to obtain
\begin{eqnarray}
(K+1) \left( \Phi_\mu(f_{K+1})-\Phi_\mu(f_\ast) \right) &\leq &  \sum_{k=0}^K \left( \Phi_\mu(f_{k+1})-\Phi_\mu(f_\ast) \right) \nonumber\\
&\leq &  \sum_{k=0}^K\frac{1}{2\tau\sigma_{\min}} \left( \|f_k-f_\ast\|^2-\|f_{k+1}-f_\ast\|^2 \right) \nonumber\\
&= & \frac{1}{2\tau\sigma_{\min}} \left( \|f_0-f_\ast\|^2-\|f_{K+1}-f_\ast\|^2 \right) \nonumber\\
&\leq &  \frac{1}{2\tau\sigma_{\min}} \|f_0-f_\ast\|^2 
\end{eqnarray}
which completes the proof.
\hfill \QED

\begin{remark}
As the existence of a converging subsequence follows from Proposition \ref{prop_bddness}, the inequality (\ref{tmp1}) with $f=f_\ast$ (the limit of a subsequence) implies $\|f_{k+1}-f_\ast\|\leq \|f_k-f_\ast\|$ and therefore convergence of $\{f_k\}_k$ under the condition $0<\tau\leq 1/\|A\|^2$ without the need to resort to the Kurdyka-Łojasiewicz inequality.
\end{remark}

\section{Extension to the group square-root Lasso}
\label{group_lasso}

The group Lasso aims at enforcing sparsity of groups of variables through a penalty that is the sum of the (quadratic) norms of 
each group \cite{Yuan2005, Giraud2015, VandeGeer2016, Stucky2017}. We assume that there is no overlap between the different groups of variables in $f$, which we denote by $f^{(j)}, j= 1, \dots, J$, the $d$-dimensional vector $f$ being the concatenation of the $f^{(j)}$'s. Hence the cost function to minimise for the group square-root Lasso \cite{Bunea2014} is given by
\begin{equation}
 \label{phi_groupsqrt}
\Phi^{\rm G}_{\mu}(f)
= \| Af-g \| + \mu\, \sum_j \Vert f^{(j)}\Vert.
\end{equation}
For simplicity, we have used the same regularisation parameter $\mu$ for each group. This is not restrictive, however, and we could easily use a different parameter $\mu_j$ for each of the groups. 

As above, we see that
\begin{equation}
\label{var_group}
{\rm arg min}_f\Phi^{\rm G}_{\mu}(f)
=  {\rm arg min}_{f ,\sigma > 0} \Phi^{\rm G}_{\mu}(f,\sigma)
\end{equation}
where
\begin{equation}
 \label{phisig_group}
\Phi^{\rm G}_{\mu}(f,\sigma)
=   \frac{\| Af-g \|^2}{\sigma} + \sigma + 2\mu\, \sum_j \Vert f^{(j)}\Vert
.\end{equation}
Starting from arbitrary  $f_0^{(j)}$'s such that the corresponding residual $\sigma_0 = \| Af_0-g \| > 0$, we can define the alternating iteration for $f$
\begin{equation}
 \label{minf2_group}
f_{k+1}
=  {\rm arg min}_{f }\left[ \| Af-g \|^2 + 2\mu \sigma_k\, \sum_j \Vert f^{(j)}\Vert \right], \quad k=0,1, \dots,
\end{equation}
and for the residual
\begin{equation}
 \label{minsig2_grp}
\sigma_{k+1}
=  \| Af_{k+1}-g \|,
\end{equation}
and form the surrogate cost function 
\begin{equation}
\label{surk_group}
{(\Phi^{\rm G}_{\mu})}^{\scriptstyle\rm{SUR}}(f, \sigma_k; f_k) =\Phi^{\rm G}_{\mu}(f, \sigma_k) + \frac{1}{\sigma_k}\left[\frac{1}{\tau}\; \|f-f_k\|^2-\|Af-Af_k\|^2\right].
\end{equation}
In analogy with (\ref{iterates}), its minimiser $f_{k+1}$ is given by replacing the soft-thresholding operator by 
a block soft-thresholding operator whose expression can be easily derived as a ``prox''  (see e.g. \cite{Combettes2005} or \cite{Beck2017}) 
\begin{equation}
\label{block_th}
S^{\rm G}_t (x) = \frac{x}{\Vert x \Vert} {\rm max} \left\{ 0, \Vert x \Vert - \frac{t}{2} \right\}
\end{equation}
to be applied on each group.
This yields the minimiser
\begin{equation} \label{iterates_grp}
f^{(j)}_{k+1}=S^{\rm G}_{2\tau\mu\sigma_k}\!\left( f^{(j)}_{k} + \tau[A^* (g-Af_{k})]^{(j)} \right)
\end{equation}
where $f^{(j)}$ denotes the vector formed by the components of the vector $f$ belonging to the group $j$.

\begin{lemma}
\label{ineq_soft_group}
The following inequality holds for the soft-thresholding operator $S^{\rm G}_t$ defined by (\ref{block_th}):
\begin{equation}
\|S^{\rm G}_t({x})-x\|\leq \sqrt{J}\ \frac{t}{2} \quad \mbox{for any} \  x\in\mathbb{R}^d
\end{equation}
where $J$ is the number of nonoverlapping groups $x^{(j)}$ in $x$.
\end{lemma}
{\em Proof:}  For each group $j$, with components $x^{(j)}$,  we have
\begin{eqnarray}
\Vert S^{\rm G}_{t}(x^{(j)})-x^{(j)}\Vert &=&
\left\Vert \frac{x^{(j)}}{\Vert x^{(j)} \Vert} \hbox{\rm max} \left\{0,  \Vert x^{(j)} \Vert - \frac{t}{2}  \right\} - x^{(j)} \right\Vert \nonumber\\
&=& \left\vert \hbox{\rm max} \left\{0,  \Vert x^{(j)} \Vert - \frac{t}{2}  \right\} - \Vert x^{(j)}  \Vert \right\vert \nonumber\\
&=&\left\{ \begin{array}{lcl}
\Vert x^{(j)} \Vert &\mbox{if}& \Vert x^{(j)}\Vert < t/2\\
t/2& \mbox{if} & \Vert x^{(j)} \Vert \geq t/2
\end{array}
\right\} \leq  \frac{t}{2} 
\end{eqnarray}
independently of $x^{(j)}$. Hence, for $J$ groups,  we have 
\begin{equation}
\|S^{\rm G}_t{x})-x\|^2 = \sum_{j=1}^J   \Vert S^{\rm G}_{t}(x^{(j)})-x^{(j)}\Vert^2 \leq  J \left(\frac{t}{2}\right)^2.
\end{equation}
Notice that for a vector in $\mathbb{R}^d$ and $d$ groups formed by single elements, we indeed recover (\ref{ineq-soft2}).
\hfill \QED

Mutatis mutandis, the proofs of the previous section can be extended to this case in a straightforward manner, replacing the soft-thresholding operator 
$S_{t}$ by $S^{\rm G}_{t}$. This leads to the analogue for the group square-root Lasso of Theorem \ref{final_cv}, namely  
\begin{theorem} \label{final_cv_group}
When $0<\tau<2/\|A\|^2$, the sequence $\{(f_k,\sigma_k)\}_k$ defined by the iteration (\ref{iterates_grp},\ref{minsig2_grp}) converges. If $\lim_k\sigma_k>0$, then $\lim_k f_k$ is a minimizer of $\Phi^{\rm G}_{\mu}$.
\end{theorem}

\begin{remark} \label{Bunea}
We must remark that the algorithm (\ref{minsig2_grp}, \ref{iterates_grp}) has been previously written in \cite{Bunea2014}. However, the proof of convergence provided there is not correct.  Indeed, it critically relies on a fixed-point theorem for non-expansive mappings due to Opial \cite{Opial1967}. This theorem had been used in \cite{Daubechies2004} to prove convergence of the ISTA algorithm for Lasso (see also there for a simplification of Opial's proof). In the case of ISTA, the iteration map $T$ is indeed non-expansive since it is the product of the non-expansive Landweber operator $Lf = f + \tau A^*(g-Af)$, with $0<\tau<2/\|A\|^2$, and of the soft-thresholding operator ${S}_t$ with a fixed threshold $t$. However, in the case of (\ref{minsig2_grp}, \ref{iterates_grp}), the threshold varies along the iteration and depends on the residual $\sigma_k=\| Af_{k}-g \|$. Hence, with $\sigma(u)=\| Au-g \|$ and $\sigma(v)=\| Av-g \|$,  $\Vert S^{\rm G}_{2\tau\mu\sigma(u)} (u) - S^{\rm G}_{2\tau\mu\sigma(v)} (v) \Vert$ cannot be bounded by $\Vert u - v\Vert$ as it is easy to see on a one-dimensional example. This unfortunately invalidates the proof of convergence given in the Appendix of  \cite{Bunea2014}.
\end{remark}

\section{Concluding remarks}
\label{conclusion}
Puzzled by the lack of scale invariance of the Lasso, we became interested in the minimisation of the cost function (\ref{phi}) instead. Then, inspired by the ISTA algorithm for the Lasso, we wrote down the iterative algorithm studied in the present paper. We discovered later that the minimisation of (\ref{phi}) had been proposed earlier under the appellation ``square-root Lasso'' and had generated a quite extensive recent literature on the subject. Motivated by such interest, we undertook to prove the convergence of our new scheme SQRT-ISTA, convergence that did not follow from the convergence results already known for ISTA because of the varying threshold, renormalised by the current residual. The potential connection with ISTA seems to have escaped to the literature on the square-root Lasso with the exception of the paper \cite{Bunea2014} which, as explained above in Remark \ref{Bunea}, contains a proof which is incorrect. 

We consider that the contribution of our paper is essentially theoretical, 
showing that the twin problems Lasso and square-root Lasso can be solved by twin first-order iterative soft-thresholding algorithms, 
enlightening the connection of our algorithm with the paradigmatic ISTA scheme and establishing its convergence under mild assumptions. The algorithm is in principle infinite-dimensional although, as shown, we had to restrict ourselves to the finite-dimensional case to finalise the convergence proof. 

Some advantages of the resulting iterative scheme are that it is simple, versatile, and easy-to-implement without resorting to sophisticated software packages. Moreover, it does not required extra conditions (except for boundedness) on the operator or matrix $A$ like those used in compressed sensing. It works for compact operators or ill-conditioned matrices, which is essential to deal with inverse problems. However, we do not pretend to have an algorithm that is competitive with respect to other numerical methods proposed in the literature for solving the square-root Lasso. The numerical experiments we did to check that our algorithm was well behaving confirmed that, as expected, its performances are similar to ISTA, in terms of computing time (with a slight extra extra cost per iteration to compute the variable thresholds) and of quality of the reconstructions, often almost identical to those provided by the Lasso for the same problem.  The regularisation parameter seems also to be less sensitive to the signal-to-noise ratio compared to ISTA, which corroborates recurrent claims in the literature as in \cite{Belloni2011,Berk2024}.

The ISTA algorithm has been the seed that fostered many subsequent algorithmic developments for the Lasso. However, it is now superseded by various accelerated and improved methods like the well-known FISTA \cite{Beck2009} which allows to get a convergence rate in $1/k^2$ instead of $1/k$ for ISTA (the $1/k$ rate being also valid for our scheme -- see Proposition \ref{proposition:rate}). We refer e.g. to \cite{Loris2009}  for an assessment of the improvement provided by several methods with respect to the original ISTA scheme. We can also refer to \cite{Mayrink2024} for a numerical comparison between different algorithms proposed for the square-root Lasso, mainly smoothed or even second-order, hence difficult to compare to our unsmoothed first-order algorithms SQRT-ISTA scheme. 

Finally, we should also mention a possible drawback of our scheme for the square-root Lasso compared to ISTA for the Lasso. 
Although SQRT-ISTA always converges, for the case of overfitting (zero residual), it does not necessarily converge to a minimum point of the cost function (\ref{phi}). On the contrary ISTA for the Lasso is always guaranteed to converge to a minimiser of (\ref{phi_Lasso}). 

\begin{figure}
\centering\includegraphics{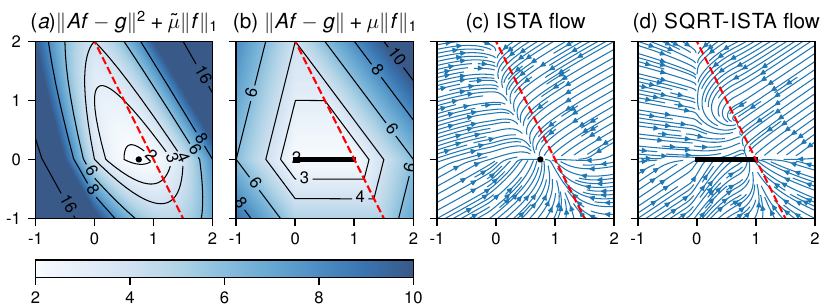}
\caption{Panels (a) and (b): Contour plots of the cost functions for the Lasso and square-root-Lasso problems for $A=(2\ 1)$ and $g=(2)$ and $\mu= 1$, $\tilde \mu = 2$. The respective minimizers are indicated with a thick  black dot or line. The set $\{f \mbox{ with } Af=g\}$ is shown with a dashed red line. Panels (c) and (d): A schematic representation of the convergence paths of the ISTA algorithm and of SQRT-ISTA algorithm (\ref{iterates_tau_summ},\ref{minsig2_summ}). While the ISTA algorithm always converges to a minimizer of the Lasso cost function, the SQRT-ISTA algorithm (\ref{iterates_tau_summ},\ref{minsig2_summ}) can also converge to a point on the red dashed line, for which $\sigma=0$.\label{figure1}}
\end{figure}

To illustrate this point, an insightful example is provided in Figure \ref{figure1} for the simple case of a 1 by 2 matrix, hence for a simple two-dimensional situation. In particular, this example illustrates the possibility of convergence of the SQRT-ISTA algorithm to a point that is not a minimiser of the cost function (\ref{phi}).

Let us observe, however, that such example is not really representative of generic realistic and properly regularised inverse imaging problems, for which the residuals are not expected to converge to zero. 

In conclusion, we believe that our work provides a valuable contribution from a theoretical point of view. From a practical perspective and focusing on inverse problems, however, we think that the square-root Lasso may present some advantages but also some drawbacks with respect to Lasso. Reaching more clearcut conclusions would require further work on the subject in relation with specific practical applications in different fields. This is beyond the scope of the present paper. 

\section*{Dedication} The present paper is dedicated to the memory of Pierre Célestin Sabatier who played a major role in establishing the field of inverse problems as a recognised independent scientific domain. Besides, as the seed motivating this work, we acknowledge a discussion on the lack of scale invariance of the Lasso with Mario Bertero, to whom we also want to bear homage.

\section*{Acknowledgments}

This work was realised with the financial support of the Fund for Scientific Research (F.R.S.–FNRS) through grant TELEVIE 40034357 (7.4586.25).

\section*{References}
\bibliography{bibsrl}
\bibliographystyle{plain}

\end{document}